\documentclass[11pt]{article}
\usepackage{color,epsfig}
\usepackage{mathrsfs}

\def\be{\begin{equation}}
\def\ee{\end{equation}}
\def\bea{\begin{eqnarray}}
\def\eea{\end{eqnarray}}
\def\bes{\begin{eqnarray*}}
\def\ees{\end{eqnarray*}}

\def\<{\langle}
\def\>{\rangle}

\def\pt{\partial}

\def\R{{\bf R}}
\def\C{{\bf C}}
\def\Z{{\bf Z}}
\def\N{{\bf N}}
\def\U{{\bf U}}

\def\bb{{\beta}}
\def\ga{{\gamma}}

\def\th{{\theta}}
\def\om{{\omega}}

\def\lm{{\lambda}}

\def\P{{\cal P}}

\def\span{{\rm span}}
\def\Sp{{\rm Sp}}

\def\hb{\vrule height0.18cm width0.14cm $\,$}
\def\ol#1{\overline{#1}}

\newtheorem{lemma}{Lemma}[section]
\newtheorem{theorem}[lemma]{Theorem}
\newtheorem{corollary}[lemma]{Corollary}

\newtheorem{proposition}[lemma]{Proposition}
\newtheorem{definition}[lemma]{Definition}

\newtheorem{remark}[lemma]{Remark}

\title{\bf Trace estimation of a family of periodic Sturm-Liouville operators with application to Robe's restricted three-body problem}

\author{Qinglong Zhou$^{1,2} $\thanks{Partially supported by the Natural Science Foundation of Zhejiang Province (No.Y19A010072) and the Fundamental Research Funds for the Central Universities (No.2017QNA3002). E-mail: \break zhouqinglong@zju.edu.cn}\quad\\
$^{1}$ Department of Mathematics\\Zhejiang University, Hangzhou 310027, Zhejiang, China\\
$^{2}$ IMCCE, Observatoire de Paris\\Avenue Denfert-Rochereau, Paris 75014, France\\
}
\date{}

\begin{document}

\maketitle

\begin{abstract}
{In this paper, we consider a family of Sturm-Liouville operators on the $\omega$-periodic domain.
The bifurcation with respect to the parameter region is studied, and the elliptic regions are estimated by trace formula.
At last, these results are used to study the linear stability of the elliptic equilibrium point along $z$-axis
in Robe's restricted three-body problem.}
\end{abstract}



\renewcommand{\theequation}{\thesection.\arabic{equation}}
\renewcommand{\thefigure}{\arabic{figure}}
\section{Introduction}

We will consider a family of periodic Sturm-Liouville operators
\begin{equation}\label{operator}
A(\beta,e)=-{d^2\over dt^2}-1+{\beta\over 1+e\cos t},
\end{equation}
where $\beta,e$ are parameters and $e\in[0,1)$, defined on the following $\omega$-periodic domain
\begin{equation}\label{om.bc.domain}
\overline{D}(\om,2\pi)=\{y\in W^{2,2}([0,2\pi],\C)\;|\;y(2\pi)=\om y(0),\dot{y}(2\pi)=\om\dot{y}(0)\}.
\end{equation}
Then it is self-adjoint.

Such an operator was studied in \cite{HLS}, \cite{Ou} and \cite{ZL1}.
We conclude their main results here:\\
(1) If $\beta>1$, $A(\beta,e)$ is positive definite on $\overline{D}(\om,2\pi)$;\\
(2) $A(1,e)$ is positive definite on $\overline{D}(\om,2\pi)$ when $\om\ne1$;\\
(3) $A(1,e)$ is semi-positive definite on $\overline{D}(1,2\pi)$ and $\ker{A(1,e)}={\rm span}\{1+e\cos t\}$.

One of the reason why to study $A(\beta,e)$ is that it has strong relation to the linear stability of periodic orbits in $N$-body problem.
For example, in three-body problem, the linear stability of Lagrangian solutions was completely determined by
the following operator ((2.29) in \cite{HLS}):
\begin{equation}\label{operator.of.Lagrangian}
A_L=-{d^2\over dt^2}I_2-I_2-{1\over2(1+e\cos t)}\left[3I_2+\sqrt{9-\beta}S(t)\right],
\end{equation}
where $S(t)=\left(\matrix{\sin2t& \cos2t\cr \cos2t& -\sin2t}\right)$,
and $\beta,e$ are parameters;
the linear stability of Euler solutions was completely determined by
the following operator ((2.42) in \cite{ZL1}):
\begin{equation}\label{operator.of.Euler}
A_E=-{d^2\over dt^2}I_2-I_2-{1\over2(1+e\cos t)}\Big[(\beta+3)I_2+3(1+\beta)S(t)\Big],
\end{equation}
where $S(t)$ is the same matrix as above, and $\beta,e$ are parameters.

However, $A_L$ and $A_E$ are more complected than $A(\beta,e)$.
Moreover, when $\beta$ takes some special values,
$A_L$ and $A_E$ are the direct sum of two $A(\beta,e)$s with some proper parameter $\beta$.
So it is worth to study $A(\beta,e)$ in details.
Recall that if $\beta>1$, $A(\beta,e)$ is positive definite on $\overline{D}(\om,2\pi)$,
so we will restrict the parameter $\beta$ to the segment $[0,1]$.

If $x\in\ker{A(\beta,e)}$, we have
\begin{equation}\label{2nd.LHS}
\ddot{x}(t)=-x(t)+{\beta\over1+e\cos t}x(t).
\end{equation}
It is a second order linear Hamiltonian system.
Letting $z=(\dot{x},x)^T$, we obtain
\begin{equation}
\dot{z}=J_2B(t)z,
\end{equation}
with
\begin{equation}
B(t)=B_{\beta,e}(t)
    =\left(\matrix{1& 0\cr 0& 1-{\beta\over1+e\cos t}}\right).
\end{equation}

Now let $\ga=\ga_{\beta,e}(t)$ be the fundamental solution of the following system:
\begin{eqnarray}
\dot{\ga}(t)&=&JB(t)\ga(t), \label{LHS}
\\
\ga(0)&=&I_2,
\end{eqnarray}
where $\beta\in[0,1]$ and $e\in[0,1)$.
Letting
\begin{equation}
L(t,x,\dot{x})={1\over2}|\dot{x}|^2+{1\over2}\left({\beta\over1+e\cos t}-1\right)|x|^2,\quad\forall x\in W^{1,2}(\R\slash2\pi\Z,\R),
\end{equation}
then zero is a solution of  the corresponding Euler-Lagrangian system.
By Legendrian transformation,
the corresponding Hamiltonian function is
\begin{equation}
H(t,z)={1\over2}B(t)z\cdot z,\qquad\forall z\in\R^2.
\end{equation}

Following \cite{Lon2} and \cite{Lon4},
denote by $\Sp(2)$ the symplectic group of real $2\times2$ matrices.
For any $\omega\in{\bf U}=\{z\in{\bf C}\;|\;|z|=1\}$ we can define a real function
$D_\omega(M)=\overline{\omega}det(M-\omega I_{2})$ for any $M\in\Sp(2)$.
Then we define $\Sp(2)_{\omega}^0 = \{M\in\Sp(2)\,|\, D_{\omega}(M)=0\}$ and
$\Sp(2)_{\omega}^{\ast} = \Sp(2)\setminus \Sp(2)_{\omega}^0$. The orientation of $\Sp(2)_{\omega}^0$ at any of its point
$M$ is defined to be the positive direction $\frac{d}{dt}Me^{t J}|_{t=0}$ of the path $Me^{t J}$ with $t>0$ small
enough. Let $\nu_{\omega}(M)=\dim_{\bf C}\ker_{\bf C}(M-\omega I_{2})$. Let
$\mathcal{P}_{2\pi}(2) = \{\gamma\in C([0,2\pi],\Sp(2))\;|\;\gamma(0)=I\}$ and
$\xi(t)={\rm diag}(2-\frac{t}{2\pi}, (2-\frac{t}{2\pi})^{-1})$ for $0\le t\le 2\pi$.

As in \cite{Lon4}, for $\lambda\in{\bf R}\setminus\{0\}$, $a\in{\bf R}$, $\theta\in (0,\pi)\cup (\pi,2\pi)$,
we denote respectively some normal forms of $2\times2$ symplectic matrices by
\begin{eqnarray}
    D(\lambda)=\left(\matrix{\lambda & 0\cr
                         0  & \lambda^{-1}}\right), \qquad
         R(\theta)=\left(\matrix{\cos\theta & -\sin\theta\cr
                        \sin\theta  & \cos\theta}\right),    \qquad
    N_1(\lambda, a)=\left(\matrix{\lambda & a\cr
                             0   & \lambda}\right). \nonumber
\end{eqnarray}

For every $M\in\Sp(2)$ and $\omega\in{\bf U}$, as in Definition 1.8.5 on p.38 of \cite{Lon4}, we define the
{\it $\omega$-homotopy set} $\Omega_{\omega}(M)$ of $M$ in $\Sp(2)$ by
$$  \Omega_{\omega}(M)=\{N\in\Sp(2)\,|\, \nu_{\omega}(N)=\nu_{\omega}(M)\},  $$
and the {\it homotopy set} $\Omega(M)$ of $M$ in $\Sp(2)$ by
\begin{eqnarray}
    \Omega(M)=\{N\in\Sp(2)\,&|&\,\sigma(N)\cap{\bf U}=\sigma(M)\cap{\bf U},\,{\it and}\; \nonumber\\
     &&\qquad \nu_{\lambda}(N)=\nu_{\lambda}(M)\qquad\forall\,\lambda\in\sigma(M)\cap{\bf U}\}.  \nonumber
\end{eqnarray}
We denote by $\Omega^0(M)$ (or $\Omega^0_{\omega}(M)$) the path connected component of $\Omega(M)$ ($\Omega_{\omega}(M)$)
which contains $M$, and call it the {\it homotopy component} (or $\omega$-{\it homtopy component}) of $M$ in
$\Sp(2)$. Following Definition 5.0.1 on p.111 of \cite{Lon4}, for $\omega\in{\bf U}$ and $\gamma_i\in \P_{\tau}(2)$
with $i=0, 1$, we write $\gamma_0\sim_{\omega}\gamma_1$ if $\gamma_0$ is homotopic to $\gamma_1$ via
a homotopy map $h\in C([0,1]\times [0,\tau], \Sp(2))$ such that $h(0)=\gamma_0$, $h(1)=\gamma_1$, $h(s)(0)=I$,
and $h(s)(\tau)\in \Omega_{\omega}^0(\gamma_0(\tau))$ for all $s\in [0,1]$. We write also $\gamma_0\sim \gamma_1$, if
$h(s)(\tau)\in \Omega^0(\gamma_0(\tau))$ for all $s\in [0,1]$ is further satisfied. We write $M\approx N$, if
$N\in \Omega^0(M)$.

For any $\gamma\in \mathcal{P}_{2\pi}(2)$ we define $\nu_\omega(\gamma)=\nu_\omega(\gamma(2\pi))$ and
$$  i_\omega(\gamma)=[\Sp(2)_\omega^0:\gamma\ast\xi], \qquad {\rm if}\;\;\gamma(2\pi)\not\in \Sp(2)_{\omega}^0,  $$
i.e., the usual homotopy intersection number, and the orientation of the joint path $\gamma\ast\xi$ is
its positive time direction under homotopy with fixed end points. When $\gamma(2\pi)\in \Sp(2)_{\omega}^0$,
we define $i_{\omega}(\gamma)$ be the index of the left rotation perturbation path $\gamma_{-\epsilon}$ with $\epsilon>0$
small enough (cf. Def. 5.4.2 on p.129 of \cite{Lon4}). The pair
$(i_{\omega}(\gamma), \nu_{\omega}(\gamma)) \in {\bf Z}\times \{0,1,2\}$ is called the index function of $\gamma$
at $\omega$. When $\nu_{\omega}(\gamma)=0$ ($\nu_{\omega}(\gamma)>0$), the path $\gamma$ is called
$\omega$-{\it non-degenerate} ($\omega$-{\it degenerate}).
For more details we refer to the Appendix or \cite{Lon4}.

Based on the above notation,
we will give some statements on the stability and instability of the periodic solutions
of the Hamiltonian systems via indices of the orbits.
Recall that, for $M\in\Sp(2)$, it is {\it linearly stable} if $\|M^j\|\leq C$ for some constant $C$ and all $j\in{\bf N}$.
Note that this implies $M$ is diagonalizable and the eigenvalues of $M$ are
all on the unit circle ${\bf U}$ of the complex plane.
We call $M$ to be {\it spectrally stable} if all its eigenvalues are on the unit circle.

\begin{definition}
Given a $T$-periodic solution $z(t)$ to a first order Hamiltonian system with fundamental solution $\gamma(t)$,
we say $z$ is spectrally stable (linearly stable) if $\gamma(T)$ is spectrally stable (linearly stable, respectively).
\end{definition}

In the literature there are many papers concerning the stability of the periodic solutions
of the Hamiltonian system using the Maslov-type index \cite{Eke,Lon2,Lon4}. The complete
iteration formula developed by Y.~Long and his collaborators is a very effective tool for
this purpose.

Precisely, if $M\in\Sp(2)$ is linearly stable as defined,
then there exists $P\in\Sp(2)$, such that \cite{BTZ} (p. 223, Remark(c))
\begin{equation}\label{1.12}
    P^{-1}MP=I_2\;\;{\rm or}\;\;R(\theta),
\end{equation}
where $\theta\in(0,2\pi)$.
Moreover, $\det(R(\theta)-I_2)>0$, which means that $\det(e^{-\epsilon J}P^{-1}MP-I_{2})>0$ with real $\epsilon>0$ small enough.
Additionally, if $M$ is $1$-nondegenerate, the second case of (\ref{1.12}) must hold;
if $M$ is $\omega$-nondegenerate for some $\omega=e^{\sqrt{-1}\theta_0}\in{\bf U}\backslash\{1\}$,
then $\theta\ne{\pm}\theta_0 ({\rm mod}\;2\pi)$ must hold.

The following two theorems describe the main results proved in this paper.

\begin{figure}\label{bifurcation_diagram}
\begin{center}
\includegraphics[height=7.5cm]{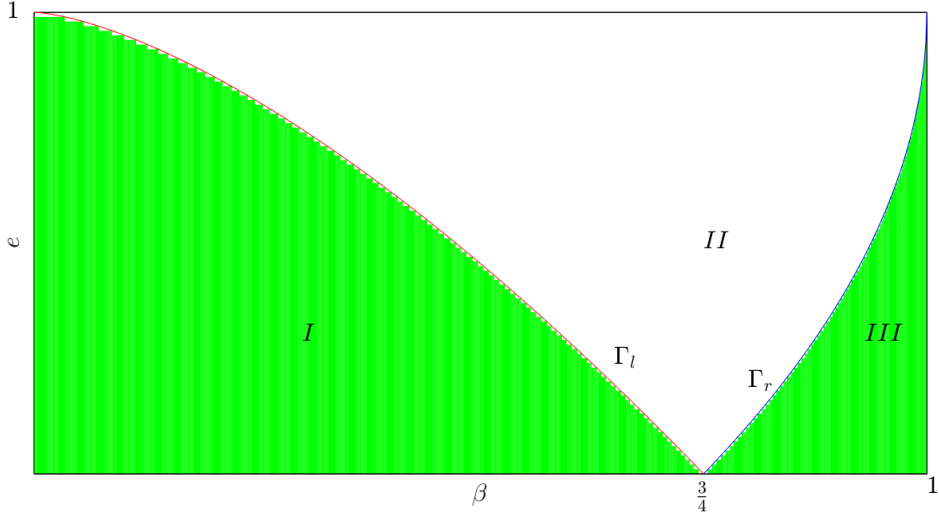}
\vskip -5 mm
\caption{The separation curves of $\ga_{\beta,e}(2\pi)$
in the $(\beta,e)$ rectangle $[0,1]\times [0,1)$}
\end{center}
\end{figure}

\begin{theorem}\label{main.theorem}
For the family of the operators $A(\beta,e)$,
denote $\ga_{\beta,e}$ by the fundamental solutions of its related first order linear Hamiltonian system.
For every $e\in[0,1)$ the $-1$-index $i_{-1}(\ga_{\beta,e})$ is non-decreasing,
and strictly decreasing only on two values $\beta=\beta_1(e)$ and $\beta=\beta_2(e)\in(0,1)$.
Define
\begin{equation}
\beta_l(e)=\min\{\beta_1(e),\beta_2(e)\}\quad{\rm and}\quad
\beta_r(e)=\max\{\beta_1(e),\beta_2(e)\}\quad{\rm for}\;e\in[0,1),
\end{equation}
and
\begin{equation}
    \Gamma_j=\{(\beta_j(e),e)\in[0,1]\times[0,1)\},
\end{equation}
for $j=l,r$.
i.e., the curves $\Gamma_l$ and $\Gamma_r$ are the diagrams of
the functions $\beta_l$ and $\beta_r$ with respect to $e\in[0,1)$, respectively.
These two curves separated the parameter rectangle $\Theta=[0,1]\times[0,1)$ into three regions,
and we denote them by I, II and III (see Figure \ref{bifurcation_diagram}), respectively.
Then we have the following:

(i) $0<\beta_i(e)<1,i=1,2$.
Moreover, $\beta_1(0)=\beta_2(0)={3\over4}$;

(ii)
The two functions $\beta_1$ and $\beta_2$ are real analytic in $e$,
and with derivatives $-{1\over2}$,
${1\over2}$ with respect to $e$ respectively,
thus they are different and the intersection points of their diagrams must be isolated if there exist when $e\in(0,1)$.
Consequently, $\Gamma_l$ and $\Gamma_r$ are different piecewise real analytic curves;

(iii) We have
\begin{equation}
    i_{-1}(\gamma_{\beta,e}) = \left\{\matrix{
   	2,&  {\it if}\;\;\beta\in[0,\beta_l(e)),\quad\;\; \cr
   	1,&  {\it if}\;\;\beta\in[\beta_l(e),\beta_r(e)),\cr
   	0,&  {\it if}\;\;\beta\in[\beta_r(e),1],\quad\;\; }\right.
\end{equation}
and $\Gamma_l,\Gamma_r$ are precisely the $-1$-degenerate curves of the path $\gamma_{\beta,e}$
in the $(\beta,e)$ rectangle $\Theta=[0,1]\times[0,1)$;

(iv) In Region I, i.e., when $0<\beta<\beta_l(e)$,
we have $\gamma_{\beta,e}(2\pi)\approx R(\theta)$ for some $\theta\in(\pi,2\pi)$,
and thus it is strongly linear stable;

(v) In Region II, i.e., when $\beta_l(e)<\beta<\beta_r(e)$,
we have $\gamma_{\beta,e}(2\pi)\approx D(-2)$,
and thus it is strongly linearly unstable;

(vi) In Region III, i.e., when $\beta_r(e)<\beta<1$,
we have $\gamma_{\beta,e}(2\pi)\approx R(\theta)$ for some $\theta\in(0,\pi)$,
and thus it is strongly linear stable.
\end{theorem}

\vspace{2mm}

\begin{remark}
For $(\beta,e)$ located on these two special curves $\Gamma_l$ and $\Gamma_r$, we have the following:

\indent(i) If $\beta_l(e)<\beta_r(e)$,
we have $\gamma_{\beta_l(e),e}(2\pi)\approx N_1(-1,1)$
Consequently, the matrix $\gamma_{\beta_l(e),e}(2\pi)$ is spectrally stable and linearly unstable;

\indent(ii) If $\beta_l(e)<\beta_r(e)$,
we have $\gamma_{\beta_r(e),e}(2\pi)\approx N_1(-1,-1)$
and it is spectrally stable and linearly unstable;

\indent(iii) If $\beta_l(e)=\beta_r(e)$,
we have $\gamma_{\beta_l(e),e}(2\pi)\approx -I_2$
and it is linearly stable.
\end{remark}

However, Figure \ref{bifurcation_diagram} was drew numerically.
Though we can prove that the bifurcation curves $\Gamma_l$ and $\Gamma_r$ are piecewise real analytic,
it is difficult to know the concrete shapes of the separation curves.
Even more, we don't know whether these two curves intersect each other at some $\beta\in(0,1)$ in the analytical sense.
In \cite{HOW}, Hu, Ou and Wang use the trace formula
to estimate the stability region and hyperbolic region for the elliptic Lagrangian orbits.
Motivated by their results, we have

\begin{theorem}\label{Elliptic_region_estimation}
$\ga_{\beta,e}$ is linearly stable if
\begin{equation}
0\le e<{1\over\sqrt{f(\beta,-1)}},\quad\beta\in[0,{3\over4}),
\end{equation}
or
\begin{equation}
0\le e<{1\over1+\sqrt{f(\beta,-1)}},\quad\beta\in({3\over4},1),
\end{equation}
where $f(\beta,\omega)$ is a function on $[0,1]\times\U$ given by (\ref{trace.function}).
\end{theorem}
In Figure 2, we have draw these two curves.

\begin{figure}\label{estimation.figure}
\begin{center}
\includegraphics[height=7.5cm]{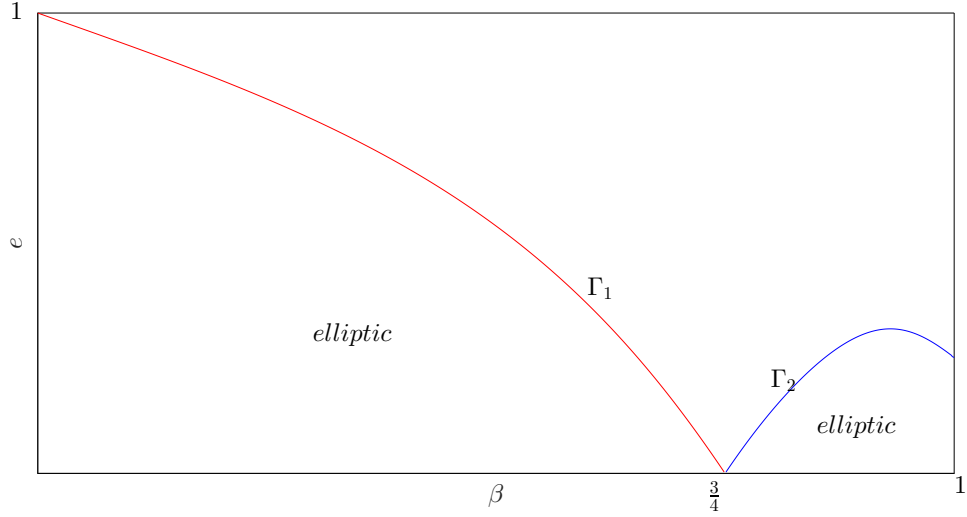}
\vskip -5 mm
\caption{The estimation of the elliptic (stable) region of $\ga_{\beta,e}(2\pi)$
given by Theorem \ref{Elliptic_region_estimation}}
\end{center}
\end{figure}
\vspace{2mm}

The paper is organized as follows.
In Section 2, we will prove Theorem \ref{main.theorem}.
In Section 3, we will prove Theorem \ref{Elliptic_region_estimation}.
In Section 4, we use the above results to study the linear stability of the elliptic equilibrium point along $z$-axis
in Robe's restricted three-body problem.

\renewcommand{\theequation}{\thesection.\arabic{equation}}
\section{Bifurcation with respect to the parameters}
\setcounter{equation}{0}

We first give the relation for the Morse index
and the Maslov-type index which covers the applications to our problem.

For $T>0$, suppose $x$ is a critical point of the functional
$$ F(x)=\int_0^TL(t,x,\dot{x})dt,  \qquad \forall\,\, x\in W^{1,2}(\R/T\Z,\R), $$
where $L\in C^2((\R/T\Z)\times \R^{2},\R)$ and satisfies the
Legendrian convexity condition $L_{p,p}(t,x,p)>0$. It is well known
that $x$ satisfies the corresponding Euler-Lagrangian
equation:
\begin{eqnarray}
    && \frac{d}{dt}L_p(t,x,\dot{x})-L_x(t,x,\dot{x})=0,    \label{A.7}\\
    && x(0)=x(T),  \qquad \dot{x}(0)=\dot{x}(T).    \label{A.8}
\end{eqnarray}

For such an extremal loop, define
\begin{eqnarray}
P(t) &=& L_{p,p}(t,x(t),\dot{x}(t)),  \nonumber\\
Q(t) &=& L_{x,p}(t,x(t),\dot{x}(t)),  \nonumber\\
R(t) &=& L_{x,x}(t,x(t),\dot{x}(t)).  \nonumber
\end{eqnarray}

For $\omega\in\U$, set
\begin{equation}
D(\omega,T)=\{y\in W^{1,2}([0,T],\C)\,|\, y(T)=\omega y(0) \}.   \label{2.10}
\end{equation}
and
\begin{equation}    \label{2.11}
    \overline{D}(\omega,T)= \{y\in W^{2,2}([0,T],\C)\,|\, y(T)=\omega y(0), \dot{y}(T)=\omega\dot{y}(0) \}.
\end{equation}

Suppose $x$ is an extreme of $F$ in $\overline{D}(\omega,T)$.
The index form of $x$ is given by
\begin{equation}
    I(y_1,y_2)=\int_0^T\{(P\dot{y}_1+Q{y_2})\cdot\dot{y}_2+Q^T\dot{y}_1\cdot{y_2}+R{y_1}\cdot{y_2}\},\;\;y_1,y_2\in{D}(\omega,T).
\end{equation}
The Hessian of $F$ at $x$ is given by
\begin{equation}
    I(y_1,y_2)=\< F''(x)y_1,y_2\>,\quad y_1,y_2\in D(\omega,T),
\end{equation}
where $\<\cdot,\cdot\>$ is the inner product in $L^2$.
Linearization of (\ref{A.7}) at $x$ is given by
\begin{equation}
-{d\over dt}(P(t)\dot{y}+Q(t)y)+Q^T\dot{y}+R(t)y=0,  \label{A.24}
\end{equation}
and $y$ is solution of (\ref{A.24}) if and only if $y\in\ker(I)$.

We define the $\omega$-Morse index $\phi_\omega(x)$ of $x$ to be the dimension of the
largest negative definite subspace of the index form $I$ which was defined on $D(\omega,T)\times D(\omega,T)$.
Moreover, $F''(x)$ is a self-adjoint operator on $L^2([0,T],\R)$ with domain $\overline{D}(\omega,T)$.
We also define
$$\nu_\omega(x)=\dim\ker(F''(x)).$$

In general, for a self-adjoint operator $A$ on the Hilbert space $\mathscr{H}$, we set
$\nu(A)=\dim\ker(A)$ and denote by $\phi(A)$ its Morse index which is the maximum dimension
of the negative definite subspace of the symmetric form $\< A\cdot,\cdot\>$. Note
that the Morse index of $A$  is equal to the total multiplicity of the negative eigenvalues
of $A$.

On the other hand, $\tilde{x}(t)=(\partial L/\partial\dot{x}(t),x(t))^T$ is the solution of the
corresponding Hamiltonian system of (\ref{A.7})-(\ref{A.8}), and its fundamental solution
$\gamma(t)$ is given by
\begin{eqnarray}
    \dot{\gamma}(t) &=& JB(t)\gamma(t),  \label{2.12}\\
    \gamma(0) &=& I_{2},  \label{2.13}
\end{eqnarray}
with
\begin{equation}
    B(t)=\left(\matrix{P^{-1}(t)& -P^{-1}(t)Q(t)\cr
                       -Q(t)^TP^{-1}(t)& Q(t)^TP^{-1}(t)Q(t)-R(t)}\right). \label{2.14}
\end{equation}

\begin{lemma}(Y.~Long, \cite{Lon4}, p.172)\label{LA.3}
For the $\omega$-Morse index $\phi_\omega(x)$ and nullity $\nu_\omega(x)$ of the solution $x=x(t)$
and the $\omega$-Maslov-type index $i_\omega(\gamma)$ and nullity $\nu_\omega(\gamma)$ of the symplectic
path $\gamma$ corresponding to $\tilde{x}$, for any $\omega\in\U$ we have
\begin{equation}
\phi_\omega(x) = i_\omega(\gamma), \qquad \nu_\omega(x) = \nu_\omega(\gamma).  \label{2.15}
\end{equation}
\end{lemma}

A generalization of the above lemma to arbitrary  boundary conditions is given in \cite{HS1}.
For more information on these topics, we refer to \cite{Lon4}.
In particular, we have for any $\bb$ and $e$,
the Morse index $\phi_{\om}(A(\bb,e))$ and nullity $\nu_{\om}(A(\bb,e))$
of the operator $A(\bb,e)$ on the domain $\ol{D}(\omega,2\pi)$ satisfy
\begin{equation}
\phi_{\om}(A(\bb,e)) = i_{\om}(\ga_{\bb,e}), \quad
\nu_{\om}(A(\bb,e)) = \nu_{\om}(\ga_{\bb,e}), \qquad
           \forall \,\om\in\U.
\end{equation}

Now we can compute the indices on the boundary segments of the parameter region.

When $\beta=0$, we have
\begin{equation}
A(0,e)=-{d^2\over dt^2}-1.
\end{equation}
Then, we obtain
\begin{equation}\label{index.of.0e}
\left\{\matrix{
i_1(\ga_{0,e})=\phi_1(A(0,e))=1,\cr
\nu_1(\ga_{0,e})=\nu_1(A(0,e))=2,
}\right.
\end{equation}
and
\begin{equation}\label{om.index.of.0e}
\left\{\matrix{
i_\om(\ga_{0,e})=\phi_\om(A(0,e))=2,\cr
\nu_\om(\ga_{0,e})=\nu_\om(A(0,e))=0,
}\right.\quad\om\in\U\backslash\{1\}.
\end{equation}

When $\beta=1$, we have
\begin{equation}
A(1,e)=-{d^2\over dt^2}-1+{1\over1+e\cos t},
\end{equation}
and such an operator has been studied in Lemma 4.1 of \cite{ZL1}.
Hence, we have
\begin{equation}\label{index.of.1e}
\left\{\matrix{
i_1(\ga_{1,e})=\phi_1(A(1,e))=0,\cr
\nu_1(\ga_{1,e})=\nu_1(A(1,e))=1,
}\right.
\end{equation}
and
\begin{equation}\label{om.index.of.1e}
\left\{\matrix{
i_\om(\ga_{1,e})=\phi_\om(A(1,e))=0,\cr
\nu_\om(\ga_{1,e})=\nu_\om(A(1,e))=0,
}\right.\quad\om\in\U\backslash\{1\}.
\end{equation}

When $e=0$, we have
\begin{equation}
A(\beta,0)=-{d^2\over dt^2}-1+\beta,
\end{equation}
and then we can compute the $\om$-indices and $\om$-nullities with respect to the parameter $\beta$ directly.
To illustrate the change of the eigenvalues of $\ga_{\beta,0}(2\pi)$ with respect to $\beta$ in detail,
we will use another method.
Noting that when $e=0$, (\ref{LHS}) becomes an ODE system with constant coefficients.
The characteristic polynomial $\det(JB-\lambda I)$ of $JB$ is given by
\begin{equation}
\lambda^2+1-\beta=0,
\end{equation}
which has two roots:
\begin{equation}
\lambda_{\pm}=\pm\sqrt{-1}\sqrt{1-\beta}.
\end{equation}
Then we get the two characteristic multipliers of the matrix $\ga_{\beta,0}(2\pi)$
\begin{equation}
\rho_\pm(\beta)=e^{\pm2\pi\sqrt{-1}\th(\beta)},
\end{equation}
where
\begin{equation}\label{theta.of.beta}
\theta(\beta)=\sqrt{1-\beta}.
\end{equation}
In particular, we obtain the following results:

When $\beta=0$, we have $\rho_\pm(0)=1$.

When $0<\beta<{3\over4}$, in (\ref{theta.of.beta}) the angle $\theta(\beta)$ decreases strictly from $1$ to ${1\over2}$
as $\beta$ increases from $0$ to ${3\over4}$.
Therefore $\rho_+(\beta)$ runs from $1$ to $-1$ clockwise along the lower semi-unit circle in the complex plane
as $\beta$ increases from $0$ to ${3\over4}$.
Correspondingly, $\rho_-(\beta)$ runs from $1$ to $-1$ counterclockwise along the upper semi-unit circle in the complex plane
as $\beta$ increases from $0$ to ${3\over4}$.

When $\beta={3\over4}$, we have $\theta({3\over4})={1\over2}$, and then $\rho_\pm({3\over4})=-1$.

When ${3\over4}<\beta<1$, the angle $\theta(\beta)$ decreases strictly from ${1\over2}$ to $0$
as $\beta$ increases from ${3\over4}$ to $1$.
Thus $\rho_+(\beta)$ runs from $-1$ to $1$ clockwise along the upper semi-unit circle in the complex plane
as $\beta$ increases from ${3\over4}$ to $1$.
Correspondingly, $\rho_-(\beta)$ runs from $-1$ to $1$ counterclockwise along the lower semi-unit circle in the complex plane
as $\beta$ increases from ${3\over4}$ to $1$.

When $\beta=1$, we have $\rho_\pm(1)=1$.

Specially, for $\om=-1$, we get
\begin{eqnarray}
i_{-1}(\ga_{\beta,0})&=&
\left\{\matrix{
2,\emph{if}\;\beta\in[0,{3\over4}),\cr
0,\emph{if}\;\beta\in[{3\over4},1],
}\right.  \label{-1.index.of.b0}
\\
\nu_{-1}(\ga_{\beta,0})&=&
\left\{\matrix{
2,\emph{if}\;\beta={3\over4},\qquad\quad\;\cr
0,\emph{if}\;\beta\in[0,1]\backslash\{{3\over4}\}.
}\right.\label{null.-1.index.of.b0}
\end{eqnarray}

Recall that $A(\bb,e)$ has a discrete spectrum and since $(A(\bb,e))_{\bb}$ is a holomorphic family
of self-adjoint operators, every eigenvalue of $(A(\bb,e))_{\bb}$ splits into
one or several eigenvalues of $A(\bb,e)$ which are holomorphic with respect to $\bb$.
Now motivated by Lemma 4.4 in \cite{HLS} and modifying its proof to our case,
we get the following lemma:

\begin{lemma}\label{Lemma:decreasing.of.index}
(i) For each fixed $e\in [0,1)$, the operator $A(\bb,e)$ is increasing
with respect to $\beta\in [0,1]$ for any fixed $\omega\in\U$. Specially
\be  \frac{\pt}{\pt\beta}A(\beta,e)|_{\bb=\bb_0} = {1\over1+e\cos t},  \ee
is a positive definite operator.

(ii) For every eigenvalue $\lm_{\bb_0}=0$ of $A(\bb_0,e_0)$ with $\om\in\U$ for some
$(\bb_0,e_0)\in [0,1]\times [0,1)$, there holds
\be \frac{d}{d\bb}\lm_{\bb}|_{\bb=\bb_0} > 0.  \ee

(iii)  For every $(\bb,e)\in(0,1)\times[0,1)$ and $\om\in\U$,
there exist $\epsilon_0=\epsilon_0(\bb,e)>0$ small enough such that for all $\epsilon\in(0,\epsilon_0)$ there holds
\be
i_\om(\ga_{\bb-\epsilon,e})-i_\om(\ga_{\bb,e})=\nu_\om(\ga_{\bb,e}).
\ee
\end{lemma}

Consequently we arrive at
\begin{corollary}\label{C2.3}
For every fixed $e\in [0,1)$ and $\omega\in \U$, the index function
$\phi_{\omega}(A(\beta,e))$, and consequently $i_{\omega}(\gamma_{\beta,e})$, is non-increasing
as $\beta$ increases from $0$ to $1$.
When $\omega=1$, these index functions are constantly equal to $1$ as $\beta\in[0,1)$,
and when $\omega\in\U\setminus\{1\}$, they are decreasing and tends from $2$ to $0$.
\end{corollary}

{\bf Proof.}
For $0<\beta_1<\beta_2\le1$ and fixed $e\in [0,1)$, when $\beta$ increases from $\beta_1$ to
$\beta_2$, it is possible that negative eigenvalues of $A(\beta_1,e)$ pass through $0$ and become
positive ones of $A(\beta_2,e)$, but it is impossible that positive eigenvalues of
$A(\beta_1,e)$ pass through $0$ and become negative by $(ii)$ of Lemma \ref{Lemma:decreasing.of.index}.
Therefore the first claim holds.

For the second claim, if $\om=1$,
we know that $i_1(\ga_{\beta,e})\le1$ when $(\beta,e)\in[0,1)\times[0,1)$ by Lemma \ref{Lemma:decreasing.of.index} $(iii)$.
Moreover, we have $\<A(\beta,e)(1+e\cos t),1+e\cos t\>=\<-1+\beta,1+e\cos t\>=-2\pi(1-\beta)<0$
when $\beta\in[0,1)$, and hence $\phi_1(A(\beta,e))\ge1$.
Therefore, we have $i_1(\ga_{\beta,e})=\phi_1(A(\beta,e))=1$.
Moreover, we must have $\nu_1(\ga_{\beta,e})=\nu_1(A(\beta,e))=0$,
otherwise we have a contradiction with Lemma \ref{Lemma:decreasing.of.index} $(iii)$.
The other cases follows from (\ref{om.index.of.0e}) and (\ref{om.index.of.1e}).
\hb

By a similar analysis to the proof of Proposition 6.1 in \cite{HLS},
for every $e\in[0,1)$ and $\omega\in\U\backslash\{1\}$,
the total multiplicity of $\omega$-degeneracy of $\gamma_{\beta,e}(2\pi)$ for $\beta\in[0,1]$ is always precisely 2, i.e.,
\begin{equation}
    \sum_{\beta\in[0,1]}\nu_\omega(\gamma_{\beta,e}(2\pi))=2,\quad \forall \omega\in\U\backslash\{1\}.
\end{equation}

Consequently,
together with the positive definiteness of $A(1,e)$ for the $\omega\in\U\backslash\{1\}$ boundary condition,
we have
\begin{theorem}\label{Th:om.degenerate.curves}
For any $\omega\in\U\backslash\{1\}$, there exist two analytic $\omega$-degenerate curves $(\beta_i(e,\omega),e)$
in $e\in[0,1)$ with $i=1,2$.
Specially, each $\beta_i(e,\omega)$ is areal analytic function in $e\in[0,1)$,
and $0<\beta_i(e,\omega)<1$ and $\gamma_{\beta_i(e,\omega),e}(2\pi)$ is $\omega$-degenerate for $\omega\in\U\backslash\{1\}$
and $i=1,2$.
\end{theorem}

{\bf Proof.}
We prove first that $i_\omega(\gamma_{\beta,e})=0$ when $\beta$ is near $1$.
By Lemma 4.1(ii) in \cite{ZL1}, $A(1,e)$ is positive definite on $\overline{D}(\omega,2\pi)$ when $\om\ne1$.
Therefore, there exists an $\epsilon>0$ small enough, which may depends on $e$ and $\omega$, such that
$A(\beta,e)$ is also positive definite on $\overline{D}(\omega,2\pi)$ when $1-\epsilon<\beta\le1$.
Hence $\nu_\omega(\gamma_{\beta,e})=\nu_\omega(A(\beta,e))=0$ when $1-\epsilon<\beta\le1$.
Thus we have proved our claim.

Then under similar steps to those of Lemma 6.2 and Theorem 6.3 in \cite{HLS},
we can prove the theorem.
\hb

Specially, for $\omega=-1$, $e\in[0,1)$,
$\beta_i(e,-1),i=1,2$ are the two $-1$-dgenerate curves.
They are exactly the same curves as $\beta_i(e),i=1,2$ in Theorem \ref{main.theorem} which we omit ``$-1$" there.

By (\ref{null.-1.index.of.b0}), $-1$ is a double eigenvalue of the matrix $\gamma_{{3\over4},0}(2\pi)$,
then the two curves bifurcation out from $({3\over4},0)$ when $e>0$ is small enough.
Thus
$\dim\ker A({3\over4},0)=\nu_{-1}(\gamma_{{3\over4},0})=2$.
Moreover, we have
\begin{equation}
    \ker A({3\over4},0)={\rm span}\left\{\cos{t\over2},\sin{t\over2}\right\}.
\label{ker.A.of.-1}
\end{equation}

Denote by $g$ the following operator
\begin{equation}
g(z)(t)=-z(2\pi-t).
\end{equation}
Obviously, $g^2=Id$ and $g$ is unitary on $L^2([0,2\pi],\R)$
One can check directly that
\begin{equation}
A(\beta,e)g=gA(\beta,e).
\end{equation}
Recall $E=\overline{D}(-1,2\pi)$ is given by (\ref{om.bc.domain}),
and let $E_+=\ker(g+I)$,$E_-=\ker(g-I)$.
Following the studies in Section 2.2 and especially the proof in Theorem 1.1 in \cite{HS1},
the subspaces $E_+$ and $E_-$ are $A(\beta,e)$-orthogonal, and $E=E_+\bigoplus E_-$.
In fact, the subspaces $E=E_-$ and $E=E_+$ are isomorphic
to the following subspaces $E_1$ and $E_2$ respectively:
\begin{eqnarray}
E_1&=&\{z\in W^{2,2}([-\pi,\pi],\R)\;|\;z\;{\rm is\;odd}\},\\
E_2&=&\{z\in W^{2,2}([-\pi,\pi],\R)\;|\;z\;{\rm is\;even}\}.\label{E_2}
\end{eqnarray}
For $(\beta,e)\in[0,1]\times[0,1)$, restricting $A(\beta,e)$ to $E_1$ and $E_2$ respectively, we then obtain
\begin{eqnarray}
\phi_{-1}(A(\beta,e))&=&\phi_{-1}(A(\beta,e)|_{E_1})+\phi_{-1}(A(\beta,e)|_{E_2}),\\
\nu_{-1}(A(\beta,e))&=&\nu_{-1}(A(\beta,e)|_{E_1})+\nu_{-1}(A(\beta,e)|_{E_2}).
\end{eqnarray}
Similar to Proposition 7.1 in \cite{HLS}, we have

\begin{proposition}
The $\om=-1$ degeneracy curve $(\beta_i(e,-1),e)$ is precisely the
degeneracy curve of $A(\beta,e)|_{E_i}$ for $i = 1$ or $2$.
\end{proposition}

Then we have the following theorem:

\begin{theorem}\label{Th:tangent.direction}
The tangent directions of the two curves $\Gamma_l$ and $\Gamma_r$ at the same bifurcation point $({3\over4},0)$
are given by
\begin{equation}
    \beta_l'(e)|_{e=0}=-{1\over2},\quad
    \beta_r'(e)|_{e=0}={1\over2}.
\end{equation}
\end{theorem}

{\bf Proof.}
Now let $(\beta(e),e)$ be one of such curves (say, the $E_2$ degenerate curve)
which starts from ${3\over4}$ with $e\in[0,\epsilon)$ for some small $\epsilon>0$ and $x_e\in E_2$
being the corresponding eigenvector, that is
\begin{equation}
    A(\beta(e),e)x_e=0.
\end{equation}
By (\ref{ker.A.of.-1}) and (\ref{E_2}), we have
\begin{equation}
    \ker(A(\beta,e)|_{E_2})=\ker(A(\beta,e))\cap E_2=\span\{\cos{t\over2}\}.
\end{equation}
Without loss of generality, we suppose
\begin{equation}
    x_0=\cos\frac{t}{2}.  \label{x_0}
\end{equation}
There holds
\begin{equation}
    \<A(\beta(e),e)x_e,x_e\>=0.\label{Axx-1}
\end{equation}

Differentiating both side of (\ref{Axx-1}) with respect to $e$ yields
$$ \beta'(e)\<\frac{\partial}{\partial \beta}A(\beta(e),e)x_e,x_e\> + (\<\frac{\partial}{\partial e}A(\beta(e),e)x_e,x_e\>
       + 2\<A(\beta(e),e)x_e,x'_e\> = 0,  $$
where $\beta'(e)$ and $x'_e$ denote the derivatives with respect to $e$. Then evaluating both
sides at $e=0$ yields
\begin{equation}  \label{4.73}
    \beta'(0)\<\frac{\partial}{\partial \beta}A({3\over4},0)x_0,x_0\>
      + \<\frac{\partial}{\partial e}A({3\over4},0)x_0,x_0\> = 0.
\end{equation}
Then we have
\begin{eqnarray}
\left.\frac{\partial}{\partial\beta}A(\beta,e)\right|_{(\beta,e)=({3\over4},0)}
    &=& 1,  \label{4.74}\\
\left.\frac{\partial}{\partial e}A(\beta,e)\right|_{(\beta,e)=({3\over4},0)}
    &=& -\cos t.  \label{4.75}
\end{eqnarray}
Thus from (\ref{x_0}), (\ref{4.74}) and (\ref{4.75}) we have
\begin{eqnarray}
    \<\frac{\partial}{\partial\beta}A({3\over4},0)x_0,x_0\>
    = \int_0^{2\pi}\cos^2\frac{t}{2}dt
    = \pi   \label{4.78}
\end{eqnarray}
and
\begin{eqnarray}
    \<\frac{\partial}{\partial e}A({3\over4},0)x_0,x_0\>
    = \int_0^{2\pi}-\cos{t}\cos^2\frac{t}{2}dt
    =-{\pi\over2}.    \label{4.79}
\end{eqnarray}
Therefore by (\ref{4.73}) and (\ref{4.78})-(\ref{4.79}),
we obtain
\begin{equation}
    \beta'(0) = {1\over2}.  \label{4.80}
\end{equation}
The other tangent can be computed similarly.
Thus the theorem is proved.
\hb

To prove Theorem \ref{main.theorem}, we need the results of splitting numbers.
Now we give their definition for $2\times2$ symplectic matrices:

\begin{definition} (\cite{Lon2}, \cite{Lon4})\label{D2.3}
For any $M\in\Sp(2)$ and $\omega\in\U$, choosing $\tau>0$ and $\gamma\in\P_\tau(2)$ with $\gamma(\tau)=M$,
we define
\begin{equation}
    S_M^{\pm}(\omega)=\lim_{\epsilon\rightarrow0^+}\;i_{\exp(\pm\epsilon\sqrt{-1}\omega)}(\gamma)-i_\omega(\gamma).
\end{equation}
They are called the splitting numbers of $M$ at $\omega$.
\end{definition}

Splitting numbers have the following properties:

\begin{lemma}
(Y.~Long, \cite{Lon4}, pp. 191)
Splitting numbers $S^{\pm}_M(\omega)$ are well defined,
i.e., they are independent of the choice of the path $\gamma\in\mathcal{P}_T(2)$ satisfying $\gamma(T)=M$.
For $\omega\in\U$ and $M\in\Sp(2)$,
splitting numbers $S^{\pm}_N(\omega)$ are constant for all $N=P^{-1}MP$, with $P\in\Sp(2)$.
\end{lemma}

\begin{lemma}
(Y.~Long, \cite{Lon4}, pp. 198--199)
For $M\in\Sp(2)$ and $\omega\in\U$, $\theta\in(0,\pi)$, there hold
\begin{eqnarray}
    S^{\pm}_M(\omega)&=&0,\quad if\; \omega\not\in\sigma(M),\label{split.num.of.D}\\
    S^{\pm}_M(\omega)&=&S^{\mp}_M(\bar\omega),\\
    0\le S^{\pm}_M(\omega)&\le&{\rm dim}\;\ker(M-\omega I),\\
    S^{+}_M(\omega)+S^{-}_M(\omega)&\le&{\rm dim}\;\ker(M-\omega I)^{2},\;\omega\in\sigma(M),\\
    (S^{+}_{N_1(1,b)}(1),S^{-}_{N_1(1,b)}(1))&=&
	\left\{\matrix{(1,1),\;{\rm if}\;b=0,1,\cr (0,0),\;{\rm if}\;b=-1,}\right. \\
    (S^{+}_{N_1(-1,a)}(-1),S^{-}_{N_1(-1,a)}(-1))&=&
	\left\{\matrix{(1,1),\;{\rm if}\;a=-1,0,\cr (0,0),\;{\rm if}\;a=1,\quad\;\;}\right. \\
    (S^{+}_{R(\theta)}(e^{\sqrt{-1}\theta}),S^{-}_{R(\theta)}(e^{\sqrt{-1}\theta}))&=&(0,1),\label{split.num.of.R}\\
    (S^{+}_{R(2\pi-\theta)}(e^{\sqrt{-1}\theta}),S^{-}_{R(2\pi-\theta)}(e^{\sqrt{-1}\theta}))&=&(1,0).
\end{eqnarray}
\end{lemma}
From the definition and property of splitting numbers, for any $\gamma\in\mathcal{P}_T(2)$ with
$\gamma(T)=M$, we have
\begin{equation} \label{index.formula}
    i_{\omega_0}(\gamma)=i_1(\gamma)+S^+_M(1)+\sum_{\omega}(S^+_M(\omega)-S^-_M(\omega))-S^-_M(\omega_0),
\end{equation}
where the sum runs over all the eigenvalues $\omega$ of $M$ belonging to the part of $\U^+=\{{\rm Re} z\ge0|z\in\U\}$
or $\U^-=\{{\rm Re} z\le0|z\in\U\}$
strictly located between $1$ and $\omega_0$.

Now we can give the proof of Theorem \ref{main.theorem}.

{\bf Proof.}
We just need to prove (iv)-(vi).

(iv) If $0<\beta<\beta_l(e)$, then by the definition of the degenerate curves,
Lemma \ref{Lemma:decreasing.of.index} $(iii)$ and Corollary \ref{C2.3}, we have
\begin{equation}
    i_1(\gamma_{\beta,e})=1,\quad\quad \nu_1(\gamma_{\beta,e})=0, \label{1.index}
\end{equation}
and
\begin{equation}
    i_{-1}(\gamma_{\beta,e})=2,\quad\quad \nu_{-1}(\gamma_{\beta,e})=0.  \label{-1.index}
\end{equation}
Then we can suppose $\gamma_{\beta,e}(2\pi)\approx M$ where $M$ is a basic normal form in $\Sp(2)$.
Moreover, we have $M=R(\theta)$ or $M=D(\pm2)$.

If $M=D(\pm2)$, by (\ref{index.formula}) and (\ref{split.num.of.D}),
we have
\begin{equation}
i_{-1}(\gamma_{\beta,e})=i_{1}(\gamma_{\beta,e}),
\end{equation}
which contradicts (\ref{1.index}) and (\ref{-1.index}).

If $M=R(\theta)$ and $\theta\in(0,\pi)$, by (\ref{index.formula}),
we have
\begin{eqnarray}
    i_{-1}(\gamma_{\beta,e})&=&i_1(\gamma_{\beta,e})-S_{R(\theta)}^-(e^{\sqrt{-1}\theta})+S_{R(\theta)}^+(e^{\sqrt{-1}\theta}) \nonumber\\
    &=&i_1(\gamma_{\beta,e})-1+0  \nonumber\\
    &=&0,
\end{eqnarray}
which contradicts (\ref{-1.index}).
Thus we must have $M=R(\theta)$ and $\theta\in(\pi,2\pi)$.

(v) If $\beta_l(e)<\beta<\beta_r(e)$, then by the definition of the degenerate curves,
Lemma \ref{Lemma:decreasing.of.index} $(iii)$ and Corollary \ref{C2.3}, we have
\begin{equation}
    i_1(\gamma_{\beta,e})=1,\quad\quad \nu_1(\gamma_{\beta,e})=0, \label{(v)1.index}
\end{equation}
and
\begin{equation}
    i_{-1}(\gamma_{\beta,e})=1,\quad\quad \nu_{-1}(\gamma_{\beta,e})=0.  \label{(v)-1.index}
\end{equation}
Then we can suppose $\gamma_{\beta,e}(2\pi)\approx M$ with $M=R(\theta)$ or $M=D(\pm2)$.

If $M=R(\theta)$, by (\ref{index.formula}) and (\ref{split.num.of.R}),
we have
\begin{eqnarray}
    |i_{-1}(\gamma_{\beta,e})-i_1(\gamma_{\beta,e})|
    =\left|-S_{R(\theta)}^-(e^{\sqrt{-1}\theta})+S_{R(\theta)}^+(e^{\sqrt{-1}\theta})\right|
    =1,
\end{eqnarray}
which contradicts (\ref{(v)1.index}) and (\ref{(v)-1.index}).

If $M=D(2)$, by Theorem 8.1.6 in \cite{Lon4},
$i_1(\ga_{\beta,e})$ must be even, which contradicts (\ref{(v)1.index}).
Thus we must have $M=D(-2)$.

(vi) can be proved similarly as $(iv)$.
\hb

\section{Estimation by the trace formula}
\setcounter{equation}{0}

We first give a brief introduction of the trace formula for linear Hamiltonian systems.
For more details, we refer to \cite{HOW}.

Consider the eigenvalue problem of the following linear Hamiltonian system,
\begin{equation}\label{eigen.problem}
\dot{z}(t)=J(B(t)+\lambda D(t))z(t),\quad z(0)=z(2\pi).
\end{equation}
Here $B(t)$ and $D(t)$ are two symmetric matrices for any $t\in[0,2\pi]$.
Denote by $A=-J{d\over dt}$, which is defined on the domian
\begin{equation}
D=\left\{z\in W^{1,2}([0,2\pi],\C^2)|\;z(0)=z(2\pi)\right\}.
\end{equation}
Then $A$ is a self-adjoint operator with compact resolvent.
Moreover, for $\lambda\in\rho(A)$, the resolvent set of $A$, $(\lambda-A)^{-1}$ is Hilbert-Schmidt.

As above, let $\ga_\lambda(t)$ be the fundamental solution of (\ref{eigen.problem}).
To state the trace formula for the Hamiltonian system, we need some notations.
Write $\hat{D}(t)=\ga_0^T(t)D(t)\ga_0(t)$.
For $k\in\N$, let
\begin{equation}
M_k=\int_0^{2\pi}J\hat{D}(t_1)\int_0^{t_1}J\hat{D}(t_2)\ldots\int_0^{t_{k-1}}J\hat{D}(t_k)dt_k\ldots dt_2dt_1,
\end{equation}
and
\begin{equation}
\mathcal{M}(\nu)=\ga_0(2\pi)\left(\ga_0(2\pi)-e^{2\pi\nu}I_2\right)^{-1},\quad G_k(\nu)=M_k\mathcal{M}(\nu).
\end{equation}

For $\nu\in\C$ such that $-J{d\over dt}-\nu J-B$ is invertible, we set
\begin{equation}
\mathcal{F}(\nu,B,D)=D(-J{d\over dt}-\nu J-B)^{-1}.
\end{equation}
In what follows, $G_k(\nu)$ and $\mathcal{F}(\nu,B,D)$ will be written in short form as $G_k$
and $\mathcal{F}$ respectively, if there is no confusion.

However, the operator $\mathcal{F}$ comes from the following boundary value problem naturally
\begin{equation}\label{eigen.problem.om_bd}
\dot{z}(t)=J(B(t)+\lambda D(t))z(t),\quad z(0)=\om z(2\pi),
\end{equation}
where $\lambda\in\R\backslash\{0\}$ and $\om=e^{2\pi\nu}$.
In fact, if we set $A_\om=-J{d\over dt}$ with the domain with $\om$-boundary condition,
then $e^{-\nu{t}}A_\om e^{\nu t}=A-\nu{J}$.
Thus $z\in\ker(A_\om-B-\lambda D)$ if and only if $e^{-\nu{t}}z(t)\in\ker(A-\nu{J}-B-\lambda D)$,
which is equivalent to ${1\over\lambda}$ being an eigenvalue of $\mathcal{F}$
provided that $A-\nu{J}-B$ is invertible.

For $m=1$, $\mathcal{F}$ is not a trace class operator but a Hilbert-Schmidt operator,
and hence $Tr(\mathcal{F})$ is not the usual trace but a kind of conditional trace \cite{HW}.
For $m\ge2$, $\mathcal{F}^m$ are trace class operators.
Note that $\lambda$ is a non-zero eigenvalue of system (\ref{eigen.problem.om_bd})
if and only if ${1\over\lambda}$ is an eigenvalue of $\mathcal{F}$.
Therefore, if the sequence $\{\lambda_i\}$ be the set of non-zero eigenvalues of the system (\ref{eigen.problem.om_bd}),
\begin{equation}
Tr(\mathcal{F}^m)=\sum_{j}{1\over\lambda_j^m},
\end{equation}
where the sum is taken for the eigenvalue ${1\over\lambda_j}$ of $\mathcal{F}$ counting the algebraic multiplicity.

We have the trace formulas:
\begin{theorem}
(\cite{HOW} Theorem 1.1 and Corollary 1.2)
For any positive integer $m$, we have
\begin{equation}\label{trace.formula}
Tr(\mathcal{F}^m)=mTr(\mathcal{G}_m)
\end{equation}
where $\mathcal{G}_m=\sum_{k=1}^m{(-1)^k\over k}\left[\sum_{j_1+\ldots+j_k=m}(G_{j_1}\ldots G_{j_k})\right]$.

Specially, we have
\begin{equation}\label{trace.formula.m=2}
Tr(\mathcal{F}^2)=Tr\left[(M_1\mathcal{M}(\nu))^2-2M_2\mathcal{M}(\nu)\right].
\end{equation}
\end{theorem}


Now set
\begin{equation}
D_{\beta,e}(t)=B_{\beta,e}(t)-B_{\beta,0}(t)={e\cos t\over1+e\cos t}K_\beta,
\end{equation}
where $K_\beta=diag(0,\beta)$.
Then we have $-J{d\over dt}-B_{\beta,e}=-J{d\over dt}-B_{\beta,0}-D_{\beta,e}$.
Let $\cos^{\pm}(t)=(\cos t\pm|\cos t|)\slash2$, and denote
\begin{equation}
K_\beta^{\pm}=\cos^{\pm}(t)K_\beta,
\end{equation}
which can be considered as two bounded self-adjoint operators on $L^2([0,2\pi],\C^2)$.
Moreover, $K_\beta^+\ge0$ and $K_\beta^-\le0$.
It is obvious that
\begin{equation}
-J{d\over dt}-\nu J-B_{\beta,0}-{e\over1-e}K_\beta^-
\ge
-J{d\over dt}-\nu J-B_{\beta,e}
\ge
-J{d\over dt}-\nu J-B_{\beta,0}-eK_\beta^+,
\end{equation}
where $\nu$ is a pure imaginary number.
Equivalently, we have
\begin{equation}
\mathcal{A}(\beta,0,\nu)-{\beta e\over1-e}\cos^-(t)
\ge\mathcal{A}(\beta,e,\nu)
\ge\mathcal{A}(\beta,0,\nu)-\beta e\cos^+(t).
\end{equation}
Here
\begin{equation}
\mathcal{A}(\beta,e,\nu)=-({d\over dt}+\nu)^2-1+{\beta\over1+e\cos t}.
\end{equation}

\begin{lemma}\label{equality.of.traces}
For an imaginary number $\nu$, such that $-J{d\over dt}-\nu J-B_{\beta,0}$ is invertible, we have
\begin{equation}
Tr\left[\mathcal{F}(\nu,B_{\beta,0},K_\beta^+)^2\right]=Tr\left[\mathcal{F}(\nu,B_{\beta,0},K_\beta^-)^2\right].
\end{equation}
\end{lemma}

{\bf Proof.}
Define an operator $G:\;x(t)\rightarrow x(t+\pi)$ on the domain $\overline{D}(1,2\pi)$, then $G^2=Id$.
We have $K_{\beta}G=GK_{\beta}$ since $K_\beta$ is a constant matrix.
Moreover, we have
\begin{equation}
(-J{d\over dt}-\nu J-B_{\beta,0})^{-1}G=G(-J{d\over dt}-\nu J-B_{\beta,0})^{-1}
\end{equation}
Therefore,
\begin{eqnarray}
&&Tr\left[\left(G\cos^+(t)K_\beta(-J{d\over dt}-\nu J-B_{\beta,0})^{-1}G\right)^2\right]
\nonumber\\
&&\qquad=Tr\left[\left(G\cos^+(t)GK_\beta(-J{d\over dt}-\nu J-B_{\beta,0})^{-1}\right)^2\right]
\nonumber\\
&&\qquad=Tr\left[\left(\cos^-(t)K_\beta(-J{d\over dt}-\nu J-B_{\beta,0})^{-1}\right)^2\right]
\end{eqnarray}
\hb

Under the assumption of Lemma \ref{equality.of.traces}, we denote
\begin{equation}\label{trace.function}
f(\beta,\om)=Tr\left[\mathcal{F}(\nu,B_{\beta,0},K_\beta^+)^2\right]=Tr\left[\mathcal{F}(\nu,B_{\beta,0},K_\beta^-)^2\right],
\end{equation}
where $\om=e^{2\pi\nu}$, is a positive function.
Motivated by Theorem 5.3 and Theorem 5.4 in \cite{HOW},
we have

\begin{theorem}
When $\beta\in[0,{3\over4})$, $\ga_{\beta,e}$ is spectrally stable if
\begin{equation}
0\le e < {1\over1+\sqrt{f(\beta,-1)}};  \label{condition1}
\end{equation}
when $\beta\in({3\over4},1)$, $\ga_{\beta,e}$ is spectrally stable if
\begin{equation}
0\le e < f(\beta,-1)^{-{1\over2}}.  \label{condition2}
\end{equation}
\end{theorem}

{\bf Proof.}
The most difference in our case is that $\ga_{\beta,e}$ is a symplectic path in $\Sp(2)$ rather than in $\Sp(4)$.
Similarly, when $\beta\in[0,{3\over4})$ and (\ref{condition1}) holds,
we have $i_{-1}(\ga_{\beta,e})\ge\ga_{-1}(\ga_{\beta,0})=2$;
when $\beta\in({3\over4},1)$ and (\ref{condition2}) holds, we have $i_{-1}(\ga_{\beta,e})=0$.
In both cases, by (4.6) of \cite{HOW}, $e(\ga_{\beta,e})\slash2\ge|i_{-1}(\ga_{\beta,e})-i_1(\ga_{\beta,e})|\ge1$,
i.e., the total algebraic multiplicity of all eigenvalues of $\ga_{\beta,e}(2\pi)$ on $\U$ is no less than $2$,
and hence the desired results are proved.
\hb

Therefore Theorem \ref{Elliptic_region_estimation} can be proved.

Let $P_\beta$ be the following $2\times2$ symplectic matrix:
\begin{equation}
P_\beta=\left(\matrix{(1-\beta)^{1\over4}& 0\cr 0& (1-\beta)^{-{1\over4}}}\right).
\end{equation}
Then we have $P_\beta^TJP_\beta=J$.

Recall $\theta(\beta)=\sqrt{1-\beta}$. Direct computation shows that
\begin{equation}
P_\beta^{-1}JB_{\beta,0}P_\beta=JP_\beta^TB_{\beta,0}P_\beta=\theta(\beta)J,
\end{equation}
for $\beta\in[0,1]$,
and hence
\begin{equation}
P_\beta^{-1}\ga_{\beta,0}(t)P_\beta=R(\theta(\beta)t).
\end{equation}
In order to diagonalize $P_\beta^{-1}\ga_{\beta,0}(t)P_\beta$,
we have
\begin{equation}
UP_\beta^{-1}\ga_{\beta,0}(t)P_\beta U^{-1}=e^{i\Theta t},
\end{equation}
where $U={1\over\sqrt{2}}\left(\matrix{1& \sqrt{-1}\cr 1& -\sqrt{-1}}\right)$
and $\Theta=diag(\theta(\beta),-\theta(\beta))$.
Changing the basis by $P_\beta U^{-1}$,
by (\ref{trace.formula.m=2}), we have
\begin{equation}
f(\beta,\om)=Tr\left[\mathcal{F}(\nu,B_{\beta,0},K_\beta^-)^2\right]=2f_1(\beta,\om)-f_2(\beta,\om)
\end{equation}
with
\begin{eqnarray}
f_1(\beta,\om)&=&Tr[-M_2\mathcal{M}(\nu)]
\nonumber\\
&=&\int_0^{2\pi}\int_0^{t}Tr\left[\cos^-(t)e^{i\Theta(s-t)}JF_\beta
      \cdot\cos^-(s)e^{i\Theta(t-s)}JF_\beta\cdot M_\beta(\om)\right]ds\;dt,
\end{eqnarray}
and
\begin{eqnarray}
f_2(\beta,\om)&=&Tr\left[-(M_1\mathcal{M}(\nu))^2\right]
\nonumber\\
&=&\int_0^{2\pi}\int_0^{2\pi}Tr\left[\cos^-(t)e^{i\Theta(s-t)}JF_\beta M_\beta(\om)
      \cdot\cos^-(s)e^{i\Theta(t-s)}JF_\beta M_\beta(\om)\right]ds\;dt,
\end{eqnarray}
where $F_\beta=U^{-T}P_\beta^TK_\beta P_\beta U^{-1}={1\over2(1-\beta)^{1\over2}}\left(\matrix{-\beta&\beta\cr\beta&-\beta}\right)$
and $M_\beta(\om)=diag\left({e^{2\pi\sqrt{-1}\th(\beta)}\over e^{2\pi\sqrt{-1}\th(\beta)}-\om},{e^{-2\pi\sqrt{-1}\th(\beta)}\over e^{-2\pi\sqrt{-1}\th(\beta)}-\om}\right)$.
Letting $\om=-1$, we have
\begin{eqnarray}
f_1(\beta,-1)&=&\int_0^{2\pi}\int_0^{t}Tr\left[\cos^-(t)
       \left(\matrix{e^{2\pi\sqrt{-1}\th(s-t)}& 0\cr 0& e^{-2\pi\sqrt{-1}\th(s-t)}}\right)
       {1\over2\sqrt{1-\beta}}\left(\matrix{-\beta& \beta\cr -\beta& \beta}\right)\right.
       \nonumber\\
       &&\qquad\qquad\;\;\cdot\cos^-(s)\left(\matrix{e^{2\pi\sqrt{-1}\th(t-s)}& 0\cr 0& e^{-2\pi\sqrt{-1}\th(t-s)}}\right)
         {1\over2\sqrt{1-\beta}}\left(\matrix{-\beta& \beta\cr -\beta& \beta}\right)
       \nonumber\\
       &&\qquad\qquad\;\;\cdot\left.\left(\matrix{e^{2\pi\sqrt{-1}\th}\over {e^{2\pi\sqrt{-1}\th}+1}& 0\cr 0& e^{-2\pi\sqrt{-1}\th}\over {e^{-2\pi\sqrt{-1}\th}+1}}\right)\right]ds\;dt
\nonumber\\
&=&{\beta^2\over4(1-\beta)}\int_{\pi\over2}^{3\pi\over2}\int_{\pi\over2}^{t}
    \cos s\cos t\left[1-{\cos(4\pi\theta(\beta)(s-t+{1\over4}))\over\cos(\pi\th(\beta))}\right]dsdt
\end{eqnarray}
and
\begin{eqnarray}
f_2(\beta,-1)
&=&\int_0^{2\pi}\int_0^{2\pi}Tr\left[\cos^-(t)
       \left(\matrix{e^{2\pi\sqrt{-1}\th(s-t)}& 0\cr 0& e^{-2\pi\sqrt{-1}\th(s-t)}}\right)
       {1\over2\sqrt{1-\beta}}\left(\matrix{-\beta& \beta\cr -\beta& \beta}\right)\right.
       \nonumber\\
       &&\qquad\left(\matrix{e^{2\pi\sqrt{-1}\th}\over {e^{2\pi\sqrt{-1}\th}+1}& 0\cr 0& e^{-2\pi\sqrt{-1}\th}\over {e^{-2\pi\sqrt{-1}\th}+1}}\right)
       \nonumber\\
       &&\qquad\qquad\;\;\cdot\cos^-(s)\left(\matrix{e^{2\pi\sqrt{-1}\th(t-s)}& 0\cr 0& e^{-2\pi\sqrt{-1}\th(t-s)}}\right)
         {1\over2\sqrt{1-\beta}}\left(\matrix{-\beta& \beta\cr -\beta& \beta}\right)
       \nonumber\\
       &&\qquad\cdot\left.\left(\matrix{e^{2\pi\sqrt{-1}\th}\over {e^{2\pi\sqrt{-1}\th}+1}& 0\cr 0& e^{-2\pi\sqrt{-1}\th}\over {e^{-2\pi\sqrt{-1}\th}+1}}\right)\right]ds\;dt
\nonumber\\
&=&{\beta^2\over4(1-\beta)}\int_{\pi\over2}^{3\pi\over2}\int_{\pi\over2}^{3\pi\over2}
    \cos s\cos t\left[1-{\cos^2(2\pi\theta(\beta)(s-t))\over\cos^2(\pi\th(\beta))}\right]dsdt
    \nonumber\\
&=&{\beta^2\over2(1-\beta)}\int_{\pi\over2}^{3\pi\over2}\int_{\pi\over2}^{t}
    \cos s\cos t\left[1-{\cos^2(2\pi\theta(\beta)(s-t))\over\cos^2(\pi\th(\beta))}\right]dsdt,
\end{eqnarray}
which in the last equality, we used the symmetry of the integral  with respect to $s$ and $t$.
Hence
\begin{equation}
f(\beta,-1)={\beta^2\over2(1-\beta)}\int_{\pi\over2}^{3\pi\over2}\int_{\pi\over2}^{3\pi\over2}
    \cos s\cos t\left[{\cos^2(2\pi\theta(\beta)(s-t))\over\cos^2(\pi\th(\beta))}
                     -{\cos(4\pi\theta(\beta)(s-t+{1\over4}))\over\cos(\pi\th(\beta))}\right]dsdt.
\end{equation}
Then we can draw the curves $\Gamma_1$ and $\Gamma_2$ in Figure 2
with the help of Matlab.


\section{Linear stability of the elliptic equilibrium point along the $z$-axis in Robe's restricted three-body problem}
\setcounter{equation}{0}

A kind of restricted three-body problem that incorporates the effect of buoyancy
forces was introduced by Robe in \cite{Robe}.
He regarded one of the primaries as a rigid spherical shell $m_1$
filled with a homogenous incompressible fluid of density $\rho_1$.
The second primary is a mass point $m_2$ outside the shell
and the third body $m_3$ is a small solid sphere of density $\rho_3$, inside the shell,
with the assumption that the mass and radius of $m_3$ are infinitesimal.
He has shown the existence of the equilibrium point with $m_3$ at the center of the shell,
where $m_2$ describes a Keplerian orbit around it.

Furthermore, he discussed two cases of the linear stability of the equilibrium point of
such a restricted three-body problem.
In the first case, the orbit of $m_2$ around $m_1$ is circular and in the second case,
the orbit is elliptic, but the shell is empty (that is no fluid inside it)
or the densities of $m_1$ and $m_3$ are equal.
In the second case, we use ``{\it elliptic equilibrium point}" to call the equilibrium point.
In each case, the domain of stability has been investigated for the whole range
of parameters occurring in the problem.

For the elliptic case,
the studies of the linear stability of equilibrium point are much more complicated
than that of the circular case,
thus in \cite{Robe},
the bifurcation diagram of linear stability was obtained just by numerical methods.
Recently, in \cite{ZZ}, the author and Y. Zhang have studied the linear stability of the elliptic equilibrium point along the $xy$-plane
analytically by using the Maslov-type index theory.
For completeness, in this section,
we will study the linear stability of such equilibrium point along the $z$-axis analytically.

The stability along the $z$-axis is governed by Equation (29) in \cite{Robe}:
\begin{equation}\label{LHS.z_axis}
\ddot{z}+z={1-\mu\over1+e\cos t}.
\end{equation}
When $\mu=1-\beta$, (\ref{LHS.z_axis}) is exactly the second linear Hamiltonian system (\ref{2nd.LHS}),
and its corresponding second order differential operator is just $A(\beta,e)$.
Then Theorem \ref{main.theorem} can be used to study such linear stability problem.

Figure 4 of \cite{Robe} is obtained numerically.
Comparing it with our Figure \ref{bifurcation_diagram}, they are symmetric with the vertical line $\beta={1\over2}$.
It is worth noting that the estimations of the  stable regions by the trace formula are given in Theorem \ref{Elliptic_region_estimation}.

\medskip

\noindent{\bf ACKNOWLEDGMENTS}

The author is partially supported by the Natural Science Foundation of Zhejiang Province (No.Y19A010072) and the Fundamental Research Funds for the Central Universities (No.2017QNA3002).

\end{document}